\documentclass[12pt,reqno]{amsart}
\usepackage{amsfonts,amsmath,amssymb}

\usepackage{longtable}
\usepackage{csvsimple}

\usepackage{dyckpaths}

\allowdisplaybreaks%
\advance\textheight by \topskip

\usepackage{mathrsfs}

\usepackage[latin1]{inputenc}
\usepackage{graphicx}

\usepackage{tikz}
\newtheorem{theorem}{Theorem}
\usepackage{euscript}

\def\[{[\! [}
\def\]{]\! ]}

\begin{document}

\title{Deepest nodes in marked ordered trees}
\author[H.~Prodinger]{Helmut Prodinger}
\address{Helmut Prodinger, Department of Mathematical Sciences, Stellenbosch
University, 7602 Stellenbosch, South Africa, and NITheCS (National Institute for
Theoretical and Computational Sciences), South Africa}
\email{hproding@sun.ac.za}
\date{\today}

\begin{abstract}
A variation of ordered trees, where each
rightmost edge might be marked or not, if it does not lead to an endnode, is investigated.
These marked ordered trees were introduced by E. Deutsch et~al.\  to model skew Dyck paths.
We study the number of deepest nodes in such trees. Explicit generating functions
are established and the average number of deepest nodes, which approaches $\frac53$ 
when the number of nodes gets large. This is to be compared to standard ordered trees
where the average number of deepest nodes approaches $2$.

\end{abstract}


\maketitle

\section{Introduction}

In~\cite{Deutsch-italy} we find the following variation of ordered trees: Each
rightmost edge might be marked or not, if it does not lead to an endnode (leaf).
They were introduced to model skew Dyck paths using trees.

We depict a marked edge by the red colour and draw all of them of size 4 (4
nodes) in a table at the end of this introductory section.

Now we move to a symbolic equation for the marked ordered trees:
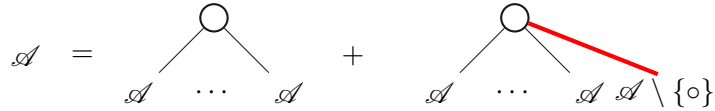
\begin{figure}[h]\small
	\begin{tikzpicture}[scale=1.0,
		s1/.style={circle=10pt,draw=black!90,thick},
		s2/.style={rectangle,draw=black!50,thick},scale=0.5]

		\node at ( 0.0,0) { $\mathscr{A}$};

		\node at (1.5,0) { $=$};

		\node(d) at (5,1)[s1]{};
		\node(e) at (3,-1){ $\mathscr{A}$};
		\node(ee) at (5,-1){$\cdots$};
		\node(f) at (7,-1){ $\mathscr{A}$};
		\path [draw,-,black!90] (d) -- (e) node{};%
		\path [draw,-,black!90] (d) -- (f) node{};%

		\node[xshift=4cm] at (0.7,0) {$+$};

		\node[xshift=4cm](d) at (5,1)[s1]{};
		\node[xshift=4cm](e) at (3,-1){ $\mathscr{A}$};
		\node[xshift=4cm](ee) at (5,-1){$\cdots$};
		\node[xshift=4cm](f) at (7,-1){ $\mathscr{A}$};
		\path [draw,-,black!90] (d) -- (e) node{};%
		\path [draw,-,black!90] (d) -- (f) node{};%
		\node[xshift=4cm](g) at (9,-1){ $\mathscr{A}\setminus\{\circ\}$};
		\path [draw,-,black,red,ultra thick] (d) -- (8.8+8,-0.5) node {};%

	\end{tikzpicture}
    \caption{Symbolic equation for marked ordered trees.\\
             $\mathscr{A}\cdots\mathscr{A}$ refers to $\ge0$ copies of
             $\mathscr{A}$.}
\end{figure}

Recall that ordered (plane, planted plane) trees are simpler and are given by
deleting the last component with the red edge.

We also bring the notion of height into the game (length of longest chain of the
root to a leaf, measured in the number of nodes). Let $\mathscr{A}_h$ denote the
family of  marked ordered trees with height $\le h$. Then
\begin{figure}[h]\small
	\begin{tikzpicture}[scale=1.0,
		s1/.style={circle=10pt,draw=black!90,thick},
		s2/.style={rectangle,draw=black!50,thick},scale=0.5]

		\node at ( 0.0,0) {$\mathscr{A}_{h+1}$};

		\node at (1.5,0) {$=$};

		\node(d) at (5,1)[s1]{};
		\node(e) at (3,-1){$\mathscr{A}_h$};
		\node(ee) at (5,-1){$\cdots$};
		\node(f) at (7,-1){$\mathscr{A}_h$};
		\path [draw,-,black!90] (d) -- (e) node{};%
		\path [draw,-,black!90] (d) -- (f) node{};%

		\node[xshift=4cm] at (0.7,0) {$+$};

		\node[xshift=4cm](d) at (5,1)[s1]{};
		\node[xshift=4cm](e) at (3,-1){$\mathscr{A}_h$};
		\node[xshift=4cm](ee) at (5,-1){$\cdots$};
		\node[xshift=4cm](f) at (7,-1){$\mathscr{A}_h$};
		\path [draw,-,black!90] (d) -- (e) node{};%
		\path [draw,-,black!90] (d) -- (f) node{};%
		\node[xshift=4cm](g) at (9,-1){$\mathscr{A}_h \setminus\{\circ\}$};
		\path [draw,-,black,red,ultra thick] (d) -- (8.8+8,-0.5) node {};%

	\end{tikzpicture}
    \caption{Symbolic equation for marked ordered trees of bounded height for
             height $h\ge2$.}
\end{figure}
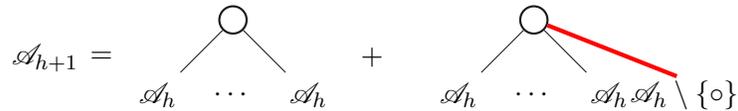

\def\bigtree{\bn{\bn{\bn{\bn{}\bn{}}\rn{\bn{}}}
             \bn{\bn{}\bn{}}\rn{\bn{}\bn{\bn{}\rn{\bn{}}}}}}
\begin{figure}[h]
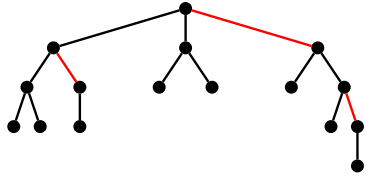
\small
    \setnodedistances{50,20,10,8,4,2,1}
    \begin{tree}\bigtree\end{tree}
    \setnodedistances{30,20,10,8,4,2,1}
    \caption{Example of a marked ordered tree.}\label{fig:bigtree}
\end{figure}
The classical bijection between ordered trees and Dyck paths consists of walking
around the tree, and recording an up-step when walking down and recording a
down-step when walking up. This can be adapted to marked ordered trees to
produce decorated Dyck paths. The additional rule is to record a red down-step
when walking up a red (marked) edge.

\begin{figure}[h]
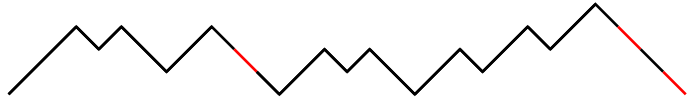
\small
    \begin{dyck}{}\bigtree\end{dyck}
    \caption{The decorated Dyck path corresponding to Figure~\ref{fig:bigtree}.}\label{fig:bigtree2}
\end{figure}

Decorated Dyck paths are in (simple) bijection to skew Dyck path, by replacing
each red down-step by a south-west (= $(-1,-1)$) step.
The next tables show all marked treed of size 4 (4 nodes) and the corresponding
objects.
\begin{figure}[h]
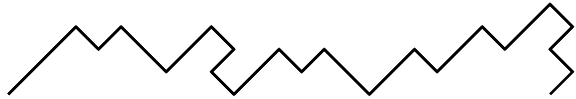
\small
    \begin{dyck}{skew}\bigtree\end{dyck}
    \caption{The skew Dyck path corresponding to Figure~\ref{fig:bigtree2}.}\label{fig:bigtree3}
\end{figure}

A representative example of trees and corresponding paths is in Figures \ref{fig:bigtree}, \ref{fig:bigtree2}, \ref{fig:bigtree3}.\footnote{Thanks are due to G. Feierabend for
the drawings.}

The main object of this paper is the analysis of the number of \emph{deepest
nodes}, i.e.\, the nodes defining the height of the tree.

For ordered trees, this was investigated by Rainer Kemp \cite{Kemp-deepest},
with important contributions provided by Volker Strehl \cite{Strehl-deepest}.

A complete list of all 10 marked ordered trees with 4 nodes is provided for the
benefit of the reader:\footnote{G. Feierabend has
compiled lists for all trees with up to 6 nodes\\ https://www.math.tugraz.at/~prodinger/pdffiles/gregg.pdf. }

\settreetablepadding{6pt}

\begin{treetable}
    \treerow{\bn{\bn{\bn{\bn{}}}}}
    \treerow{\bn{\rn{\bn{\bn{}}}}}
    \treerow{\bn{\bn{\rn{\bn{}}}}}
    \treerow{\bn{\rn{\rn{\bn{}}}}}
    \treerow{\bn{\bn{}\bn{\bn{}}}}
    \treerow{\bn{\bn{}\rn{\bn{}}}}
    \treerow{\bn{\bn{\bn{}}\bn{}}}
    \treerow{\bn{\bn{\bn{}\bn{}}}}
    \treerow{\bn{\rn{\bn{}\bn{}}}}
    \treerow{\bn{\bn{}\bn{}\bn{}}}
\end{treetable}

For completeness, it is mentioned that the average height of such marked ordered trees was already identified to be asymptotic to $\frac2{\sqrt5}\sqrt{\pi n}$ \cite{garden}, which is
slightly smaller than $\sqrt{\pi n}$ in the classical case \cite{BrKnRi72}.

\section{Enumeration}

We start by the enumerating the marked trees according to the number of nodes.
Translating the symbolic equation,
\begin{equation*}
    A = \frac{z}{1-A}+\frac{z(A-z)}{1-A}=-z+\frac{z(2-z)}{1-A},
\end{equation*}
with the relevant solution
\begin{align*}
    A(z) &= \frac{1-z-\sqrt{1-6z+5z^2}}{2}\\&
         = z+{z}^{2}+3{z}^{3}+10{z}^{4}+36{z}^{5}
           +137{z}^{6}+543{z}^{7}+2219{z}^{8}+\cdots,
\end{align*}
and the sequence $1,1,3,10,36,\dots$ of coefficients  is sequence A002212 in~\cite{OEIS}.
Next we enumerate the classes $\mathscr{A}_h$ according to the size.
The treatment of deepest nodes will come a bit later. The enumerating sequence
of $\mathscr{A}_h$ is defined to be $A_h=A_h(z)=\frac{f_h}{g_h}$. The
recursion is
\begin{equation*}
    A_{h+1}=-z+\frac{z(2-z)}{1-A_h},\quad A_1=z.
\end{equation*}
We may set  $f_1=z$, $g_1=1$, $f_2=z$, $g_2=1-z$. Then
\begin{equation*}
    f_{h+1}=zf_h+z(1-z)g_h,\quad g_{h+1}=g_h-f_h.
\end{equation*}
From this
\begin{equation*}
	g_{h+1}-g_{h+2}=zg_h-zg_{h+1}+z(1-z)g_h.
\end{equation*}
Solving the characteristic equation $X-X^2=z-zX+z(1-z)$, we find the two roots
\begin{equation*}
    \lambda=\frac{1+z+\sqrt{1-6z+5z^2}}{2},
    \quad \mu=\frac{1+z-\sqrt{1-6z+5z^2}}{2}.
\end{equation*}
The solution must be of the form
\begin{equation*}
    A_h=\frac{C_1\lambda^h-C_2\mu^h}{C_3\lambda^h-C_4\mu^h},
\end{equation*}
and an attractive form could be written using the substitution
$z=\frac{v}{1+3v+v^2}$, since then $\frac\lambda z=2+v^{-1}$ and $\frac \mu
z=2+v$.
Then
\begin{equation*}
    A_h = z(1+v)\frac{{(1+2v)}^{h-1}-v^h{(v+2)}^{h-1}}
                     {{(1+2v)}^{h-1}-v^{h+1}{(v+2)}^{h-1}},
\end{equation*}
which could be proved by induction as well.\footnote{G. Feierabend has worked
out the details of such a proof \\ https://www.math.tugraz.at/~prodinger/pdffiles/gregg.pdf.} It is also worthwhile to write
\begin{equation*}
v=\frac{1-3z-\sqrt{1-6z+5z^2}}{2z}.
\end{equation*}
Note that $\frac\mu z=2+v$
and $\frac\lambda z=2+v^{-1}$. Further
\begin{equation*}
	\frac\mu\lambda =\frac{(1+z)}{(2-z)z}\mu-1=\frac{v(2+v)}{1+2v}=:q.
\end{equation*}
All these equivalent forms are useful somehow.

Now we count the deepest nodes, using a second variable $t$. We write
$p_h=p_h(z,t)$, and $[z^{n}t^{i}]p_{h}(z,t)$ is the number of marked ordered
trees with $n$ nodes, height $\le h$, and $i$ nodes on level $h$. For $i=0$,
this means that the tree has height $<h$, so $p_h(z,0)=p_{h-1}(z,1)$. The
symbolic equation is used, but with a twist, since the recursion does not allow
to compute $p_2$. We have
\begin{equation*}
    p_1 = zt,\quad p_2=\frac{z}{1-p_1}
        = \frac{z}{1-zt},\quad p_{h+1}
        = -z+\frac{z(2-z)}{1-p_{h}},\text{ for } h\ge2.
\end{equation*}
Therefore
\begin{equation*}
    p_3=-z+\frac{z(2-z)}{1-p_2}=-z+\cfrac{z(2-z)}{1-\cfrac{z}{1-zt}},
\end{equation*}
\begin{equation*}
    p_4 = -z+\frac{z(2-z)}{1-p_3}
        = -z+\cfrac{z(2-z)}{1+z-\cfrac{z(2-z)}{1-\cfrac{z}{1-zt}}},
\end{equation*}
\begin{equation*}
    p_5=-z+\frac{z(2-z)}{1-p_4}=-z+\cfrac{z(2-z)}{1+z-
    \cfrac{z(2-z)}{1+z-\cfrac{z(2-z)}{1-\cfrac{z}{1-zt}}}}, \quad \&c.
\end{equation*}
Expanding
\begin{align*}
    p_2&=z+t z^2+t^2 z^3+t^3 z^4+\cdots,\\
    p_3&=z+z^2+(1+2 t) z^3+(1+3 t+2 t^2) z^4+\cdots,\\
    p_4&=z+z^2+3z^3+(6+4t)z^4+\cdots,\\
    p_5&=z+z^2+3z^3+10z^4+\cdots.
\end{align*}
We look at the coefficient of $z^4$ and think about the list of 10 trees
drawn earlier. For height $\le2$, one such tree appears, and it has 3
deepest nodes.
Next, 3 trees appear with one deepest node, and 2 with two deepest nodes. For
height $\le 4$, four further trees appear, with one deepest node each.

With a lot of help from \textsc{Gfun}~\cite{SaZi94}, we get for $h\ge2$
\begin{align*}
    p_h &= z(1+v)\\
        &  \times\frac{(v^2t-v^2+vt-3v-1){(1+2v)}^{h-2}
           -(-v^2+vt-3v+t-1)v^{h-1}{(2+v)}^{h-2}}{(v^2t-v^2+vt-3v-1)
           {(1+2v)}^{h-2}-(-v^2+vt-3v+t-1)v^{h}{(2+v)}^{h-2}}\\
        &= z(1+v)\frac{(v^2t-v^2+vt-3v-1)\lambda^{h-2}-(-v^2+vt-3v+t-1)
           v\mu^{h-2}}{(v^2t-v^2+vt-3v-1)\lambda^{h-2}-(-v^2+vt-3v+t-1)
           v^{2}\mu^{h-2}}\\
        &= z(1+v)\frac{1-Rvq^{h-2}}{1-Rv^{2}q^{h-2}}
\end{align*}
with
\begin{equation*}
	R = \frac{-v^2+vt-3v+t-1}{v^2t-v^2+vt-3v-1}=1
        +\frac{v-1}{v}\sum_{k\ge1}{(1+v)}^{k}t^{k}z^{k}.
\end{equation*}
The representation
\begin{equation*}
	R=1+\frac{v-1}{v}\frac{t(1+v)z}{1-t(1+v)z}
\end{equation*}
might be the most attractive.
The reader can compare this for $t=1$ with $A_h$ given earlier.

The most interesting generating function is
\begin{equation*}
G(z,t):=zt+\sum_{h\ge2}\Big(p_h(z,t)-p_h(z,0)\Big);
\end{equation*}
the coefficient of $z^{n}t^{i}$ in  $G(z,t)$ for $n,i\ge1$ is the number of
marked ordered trees with $n$ nodes and $i$ deepest nodes.

\section{Continuing with exact analysis}

First, note that $R|_{t=0}=\dfrac{-v^2-3v-1}{-v^2-3v-1}=1$. Then
\begin{align*}
    \frac{p_h(z,t)-p_h(z,0)}
    {z(1+v)}&=\frac{1-Rvq^{h-2}}{1-Rv^{2}q^{h-2}}
    -\frac{1-vq^{h-2}}{1-v^{2}q^{h-2}}\\
    &= \frac{(1-Rvq^{h-2})(1-v^{2}q^{h-2})-(1-vq^{h-2})(1-Rv^{2}q^{h-2})}
            {(1-Rv^{2}q^{h-2})(1-v^{2}q^{h-2})}\\
    &= \frac{(1-v)v(1-R)q^{h-2}}{(1-Rv^{2}q^{h-2})(1-v^{2}q^{h-2})}\\
    &= (1-v)v     \bigg[\frac{q^{h-2}}{1-v^2q^{h-2}}-
       \frac{Rq^{h-2}}{1-Rv^2q^{h-2}}\bigg],
 \end{align*}
or
\begin{align*}
    p_{h}(z,t)-p_{h}(z,0)&= \frac{z(1-v^2)}{v}
    \bigg[\frac{v^2q^{h-2}}{1-v^2q^{h-2}}-
    \frac{Rv^2q^{h-2}}{1-Rv^2q^{h-2}}\bigg].
\end{align*}
Summing,
 \begin{align*}
     \sum_{h\ge1}\Big(p_{h+1}(z,t)-p_{h+1}(z,0)\Big)&= \frac{z(1-v^2)}{v}
     \sum_{h\ge1}  \bigg[\frac{v^2q^{h-1}}{1-v^2q^{h-1}}-
     \frac{Rv^2q^{h-1}}{1-Rv^2q^{h-1}}\bigg]\\
     &= \frac{z(1-v^2)}{v}  \sum_{h\ge1}  \bigg[\frac{\delta q^{h}}{1-\delta
        q^{h}}- \frac{R\delta q^{h}}{1-R\delta q^{h}}\bigg]\\
     &= \frac{z(1-v^2)}{v}  \sum_{k\ge1}  \bigg[\frac{\delta^k q^{k}}{1- 
        q^{k}}- \frac{R^k\delta^k q^{k}}{1-  q^{k}}\bigg]\\
     &= \frac{z(1-v^2)}{v}  \sum_{k\ge1}   \frac{\delta^k q^{k}(1-R^k)}{1-
        q^{k}}
\end{align*}
 with
 \begin{equation*}
 	\delta =\frac{v(2v+1)}{v+2}.
 \end{equation*}
Using the binomial theorem,
\begin{equation*}
	R^k-1
    =\sum_{i=1}^k\binom ki{\bigg(\frac{v-1}{v}\frac{t(1+v)z}{1-t(1+v)z}\bigg)}^i.
\end{equation*}
Putting things together, 
\begin{align*}
	G(z,t)-zt &= \sum_{h\ge1}\Big(p_{h+1}(z,t)-p_{h+1}(z,0)\Big)\\*
              &= \frac{z(v^2-1)}{v}  \sum_{1\le i\le k}  \binom
                 ki{\bigg(\frac{v-1}{v}\frac{t(1+v)z}{1-t(1+v)z}\bigg)}^i
                 \frac{\delta^k q^{k}}{1-  q^{k}}.
\end{align*}
The generating function is now fully explicit.
\begin{theorem}
	The generating function $G(z,t)$ where the coefficient of $z^nt^i$ refers to the number of marked ordered trees with $n$ nodes and $i$ deepest nodes, has the explicit form
	\begin{equation*}
	G(z,t)=zt+ \frac{z(v^2-1)}{v}  \sum_{1\le i\le k}  \binom
ki{\bigg(\frac{v-1}{v}\frac{t(1+v)z}{1-t(1+v)z}\bigg)}^i
\frac{\delta^k q^{k}}{1-  q^{k}},
	\end{equation*}
with $z=\dfrac{v}{1+3v+v^2}$, $q=\dfrac{v(v+2)}{2v+1}$, and $\delta =\dfrac{v(2v+1)}{v+2}$.\qed
\end{theorem}

Now we are interested in the average number of deepest nodes, assuming all trees of size $n$ to be equally likely. For that, we have to differentiate $G(z,t)$ w.r.t. $t$, followed by
$t=1$. We ignore the tree with one node and one deepest node. Only the quantity $R$ contains the variable $t$:
\begin{align*}
    \frac{d}{dt}(1-R^h)\Big|_{t=1}
    = \frac{(1-v^2)(1+3v+v^2)h{(v+2)}^{h-1}v^{h-1}}{{(2v+1)}^{h+1}}
    = \frac{(1-v^2)(1+3v+v^2)}{(v+2)v(2v+1)}hq^h.
\end{align*}
Therefore
\begin{align*}
\frac{d}{dt}\sum_{h\ge1}\Big(p_{h+1}(z,t)-p_{h+1}(z,0)\Big)\bigg|_{t=1}&=\frac{z(1-v^2)}{v}  \sum_{k\ge1}   \frac{\delta^k q^{k}\frac{(1-v^2)(1+3v+v^2)}{(v+2)v(2v+1)}kq^k}{1-
	q^{k}}\\*
&=    \frac{{(1-v^2)}^2}{(v+2)v(2v+1)} \sum_{k\ge1}
    \frac{k\delta^k q^{2k}}{1-q^{k}}.
\end{align*}

\section{Asymptotics}

For the following, we refer to \cite{FGD} and use a hybrid approach, first the Mellin transform, to establish to local behaviour, and then singularity analysis to switch to
the behaviour of the coefficients. The book \cite{FS} is of course also relevant here.

The goal is to find the behaviour of 
\begin{equation*}
\sum_{k\ge1}\frac{k\delta^k q^{2k}}{1-q^{k}}=\sum_{k\ge1}\frac{kv^{2k} q^{k}}{1-q^{k}}
=\sum_{k,j\ge1}kv^{2k} q^{jk}
\end{equation*}
as $z\to\frac15$, or $v\to1$. First, we start with the simpler sum
\begin{equation*}
\sum_{k\ge1}\frac{k q^{3k}}{1-q^{k}}
\end{equation*}
and discuss later that the difference is negligible. We set $q=e^{-w}$.   Then we deal with
\begin{equation*}
\sum_{k\ge1}\frac{ke^{-3kw}}{1-e^{-kw}}=\sum_{k\ge1}\sum_{j\ge 3}ke^{-kjw}.
\end{equation*}
The Mellin transform of this is then $\Gamma(s)\zeta(s-1)\bigl(\zeta(s)-1-\frac1{2^s}\bigr)$, and the next step is to find the residues of
$\Gamma(s)\zeta(s-1)(\zeta(s)-1-\frac1{2^s})w^{-s}$ left to the line $\Re s=\frac32$, say. We compute these residues at $s=1$
and $s=0$ (with a computer), with the (cumulative) result $-\frac{1}{2w}+\frac{5}{24}$. But $w=-\log(q)$, and we expand
\begin{equation*}
 \frac{{(1-v^2)}^2}{(v+2)v(2v+1)}\Big(\frac{1}{2\log(q)}+\frac{5}{24}\Big)
\end{equation*}
around $v=1$, with the result
\begin{equation*}
-\frac13(1-v)-\frac{2}{27}(1-v)^2+\cdots.
\end{equation*}
 One could from this translate to an expansion about $1- 5z$, but it is not necessary, since the generating function of the
marked ordered trees, established to be  $z(1+v)$, which we need for normalization, is $\sim \frac 25-\frac15(1-v)+\cdots$, and the quotient $(-\frac 13)/(-\frac15)=\frac{5}{3}$ is the average number of deepest nodes (leading term), when all trees of size $n$ are considered to be equally likely and $n$ gets large.
\begin{theorem}
	The average number of deepest nodes, when all marked ordered trees with $n$ nodes are considered to be equally likely, approaches $\frac{5}{3}$ as $n\to\infty$.\qed
\end{theorem}

\tiny
\renewcommand{\arraystretch}{1.2}
\begin{longtable}{llll}%
	\bfseries Nodes & \bfseries Deepest nodes &
	\bfseries Trees & \bfseries Ratios
	\begin{filecontents*}{ratios.csv}
Number of nodes, Number of deepest nodes, Number of trees, Ratio
2, 1, 1, 1.00000000000000
3, 4, 3, 1.33333333333333
4, 14, 10, 1.40000000000000
5, 52, 36, 1.44444444444444
6, 202, 137, 1.47445255474453
7, 814, 543, 1.49907918968692
8, 3367, 2219, 1.51735015772871
9, 14224, 9285, 1.53193322563274
10, 61122, 39587, 1.54399171445171
11, 266336, 171369, 1.55416673960868
12, 1174054, 751236, 1.56282978983968
13, 5226196, 3328218, 1.57026853409242
14, 23459020, 14878455, 1.57671075390556
15, 106065578, 67030785, 1.58234127796653
16, 482598675, 304036170, 1.58730678326858
17, 2208111308, 1387247580, 1.59172114612736
18, 10153335117, 6363044315, 1.59567254514713
19, 46894469566, 29323149825, 1.59923029571739
20, 217453338987, 135700543190, 1.60245002617664
21, 1011990062528, 630375241380, 1.60537723580733
22, 4725078802079, 2938391049395, 1.60804968523569
23, 22127901099074, 13739779184085, 1.61049903368935
24, 103910897639245, 64430797069375, 1.61275201247876
25, 489188386162736, 302934667061301, 1.61483131299643
26, 2308345828917289, 1427763630578197, 1.61675628898215
27, 10915917653075084, 6744284275226223, 1.61854352628232
28, 51723425586415104, 31923955212096244, 1.62020730961359
29, 245539814027935212, 151403298421257630, 1.62176000515363
\end{filecontents*}
\csvreader[head to column names]{ratios.csv}{}
	{\\\hline\csvcoli&\csvcolii&\csvcoliii&\csvcoliv}
\end{longtable}

\normalsize

It remains to discuss that replacing $v$ by $q$ leads to the same main term. We have
\begin{equation*}
	v-q=\frac{v(v-1)}{1+2v},\quad v^2-q^2=\frac{3v^2(v-1)(1+v)}{(1+2v)^2},\quad v^3-q^3=\frac {{v}^{3} \left( v-1 \right)  \left( 7{v}^{2}+13v+7 \right) 
	}{ \left( 1+2v \right) ^{3}}
\end{equation*}
and so on; there is always a factor $v-1$ present, so the asymptotics of the difference of the two sums has an extra factor $v-1$, which means one order of magnitude smaller.
This is perhaps easier to see when switching to the $z$-world: $1-v\sim\sqrt5\sqrt{1-5z}$, and for instance $(1-5z)^{3/2}$ leads already to coefficients that are smaller by a factor $\frac 1n$. The necessary background information can be found in \cite{FO}.

\section{Conclusion} 

The continued fraction expression for $p_h=p_h(z,t)$ contains the variable $t$ only at the bottom level. It would be desirable to have an equivalent representation
where $t$ appears only near the top of the continued fraction. This should then lead to an identity of the Kemp/Strehl type \cite{Kemp-deepest,Strehl-deepest}.
Flajolet's paper \cite{Flajolet-aspects} does not seem to be immediately applicable.

So far, we were not successful with this.


\bibliographystyle{plain}

\begin{thebibliography}{1}

    \bibitem{BrKnRi72} N.~G. De~Bruijn, D.~E. Knuth, and S.~O. Rice.
        \newblock The average height of planted plane trees.
        \newblock In R.~C. Read, editor, {\em Graph Theory and Computing}, pages
        15--22. Academic Press, 1972.

    \bibitem{Deutsch-italy} E. Deutsch, E. Munarini, and S. Rinaldi.
        \newblock Skew Dyck paths.
        \newblock  {\em  J. Stat. Plann. Infer.}, 140 (8) (2010) 2191--2203.

    \bibitem{Flajolet-aspects} P. Flajolet.
	    \newblock Combinatorial aspects of continued fractions.
		\newblock {\em  Discrete Math.}, 32 (1980) 125--161.


\bibitem{FGD} P. Flajolet, X. Gourdon, and P. Dumas. Mellin transforms and asymptotics: harmonic sums.
Theoret. Comput. Sci. 144 (1995),  3--58.


\bibitem{FO} P. Flajolet and A. Odlyzko. Singularity Analysis of Generating Functions.
SIAM J. Discrete Math., 3 216--240. 




\bibitem{FS} P. Flajolet and R. Sedgewick. 
Analytic combinatorics. Cambridge University Press, Cambridge, 2009. 

    \bibitem{Kemp-deepest} R. Kemp.
        \newblock 	On the number of deepest nodes in ordered trees.
        \newblock  {\em  Discrete Math.}, 81 (1990) 247--258.



    \bibitem{garden} H. Prodinger.
\newblock 	A walk in my lattice path garden. preprint 2021.
\newblock  {\em https://arxiv.org/abs/2111.14797}.



    \bibitem{SaZi94} B. Salvy and P. Zimmermann, Gfun: a Maple package for the
        manipulation of generating and holonomic functions in one variable, ACM
        Transactions on Mathematical Software, vol. 20, no. 2, pp. 163--177,
        1994.

    \bibitem{OEIS} Neil J.~A. Sloane and The OEIS~Foundation Inc.
        \newblock The on-line encyclopedia of integer sequences, 2021.

    \bibitem{Strehl-deepest} V. Strehl.
        \newblock Two short proofs of Kemp's identity for rooted plane trees.
        \newblock  {\em  European J. Combin.}, 5 (1984) 373--376.

\end{thebibliography}

\end{document}